\documentclass[11pt]{amsart}

\usepackage{amsmath}
\usepackage{amssymb}
\usepackage[all]{xy}
\usepackage{sty1}
\usepackage{mathdots}

\makeatletter
\@addtoreset{equation}{section}
\makeatother

\begin{document}


\newtheorem{unique}{Proposition}[section]
\newtheorem{tpdiag}[unique]{Proposition}

\newtheorem{unit}{Proposition}[section]
\newtheorem{algorithm}[unit]{Lemma}
\newtheorem{linear}[unit]{Example}
\newtheorem{flatsem}[unit]{Example}
\newtheorem{uflatsem}[unit]{Example}
\newtheorem{powset}[unit]{Example}
\newtheorem{amenconst}[unit]{Theorem}
\newtheorem{amenconst1}[unit]{Corollary}

\newtheorem{gradedalga}{Proposition}[section]
\newtheorem{gradedalg}[gradedalga]{Theorem}
\newtheorem{gradedlinear}[gradedalga]{Proposition}
\newtheorem{dalesq}[gradedalga]{Proposition}

\title[Amenability constants]{Amenability constants for semilattice algebras}

\author{Mahya Ghandehari, Hamed Hatami and Nico Spronk}

\begin{abstract}
For any finite commutative idempotent semigroup $S$, a 
{\it semilattice}, we show how to compute the amenability
constant of its semigroup algebra $\ell^1(S)$, which is always
of the form $4n+1$.  We then
show that these give lower bounds to amenability constants
of certain Banach algebras graded over semilattices.
We also demonstrate an example of a commutative Clifford semigroup $G$
for which amenability constant of $\ell^1(G)$ is not of the form $4n+1$.
We also show there is no commutative semigroup with amenability constant
between $5$ and $9$.
\end{abstract}

\maketitle

\footnote{


Research of the third named author supported by NSERC Grant 312515-05.}


In conjunction with V.\ Runde \cite{rundes1}, the third named author proved 
that for a locally compact group $G$, $G$ is compact if and only if its
Fourier-Stieltjes algebra $\mathrm{B}(G)$ is operator amenable with operator
amenability constant less than $5$.  In a subsequent article \cite{rundes2},
examples of non-compact groups $G_1$ were found for which the operator 
amenability constant is exactly $5$.  In related work of Dales, Lau and
Strauss  \cite[Corollary 10.26]{dalesls}, improving on 
\cite[Theorem 3.2]{stokke}, it was
shown that a semigroup algebra $\ell^1(S)$ has amenability constant
less than $5$, if and only if $S$ is an amenable group.  For the multiplicative
semigroup $L_1=\{0,1\}$, it is known that the amenability constant
of $\ell^1(L_1)$ is $5$.  These parallel facts are not coincidences
since for the special groups $G_1$, mentioned above, $\mathrm{B}(G_1)$
is {\it $\ell^1$-graded} over $L_1$, i.e.\ there are $1$-operator amenable
subalgebras $\fA_0$ and $\fA_1$ such that $\mathrm{B}(G_1)=
\fA_0\oplus_{\ell^1}\fA_1$, and $\fA_0$ is an ideal.

We are thus led to consider the general situation of Banach algebras
graded over semilattices, i.e.\ commutative idempotent semigroups,
which we define in Section \ref{sec:graded}.  To do this, in Section
\ref{sec:semilattice} we
develop a method for computing the amenability constants associated
to finite semilattice algebras.  The
results in Section \ref{sec:semilattice} have a similar flavour to some results
from those in the recent monograph 
\cite{dalesls}, and are very similar to some results of Duncan and Namioka
\cite{duncann}.  However, our method is explicit and quantitative, and
thus is a nice complement to their work.  
In Section \ref{sec:graded} we obtain a lower bound for the amenability
constant of Banach algebras graded over finite semilattices.  We
show a surprising example which indicates our lower bound is not,
in general the amenability constant.  We show, at least for certain
finite dimensional algebras graded over linear semilattices, that 
our lower bound is achieved.  We close with an answer to a question
asked of us by H.G. Dales:  
we show that there does not exist a commutative semigroup
$G$ such that $5<\mathrm{AM}(\ell^1(G))<9$.

There are natural examples
of Banach algebras from harmonic analysis, due to Taylor \cite{taylor},
Inoue \cite{inoue}, and Ilie and Spronk \cite{ilies,ilies1}, to which our 
techniques 
apply.  We recommend the reader to \cite{ilies} and \cite{rundes1} for
more on this.   We feel that ideas developed here may lead to
a tool to help classify which locally compact groups admit
operator amenable Fourier-Stieltjes algebras
$\mathrm{B}(G)$.  Our hope is that the operator
amenability constants $\mathrm{AM}_{op}(\mathrm{B}(G))$
can all be computed.  We conjecture they are a subset
of $\{4n+1:n\in\En\}$, motived by Theorem \ref{theo:amenconst}
and Theorem \ref{theo:gradedalg},
below.  We hope that these values will serve as
a tool for classifying for which groups $G$, $\mathrm{B}(G)$
is operator amenable.

Interest in amenability of semigroup algebras, in particular for inverse
semigroups and Clifford semigroups, goes back at least
as far as Duncan and Namioka \cite{duncann}. 
Gr{\o}nb{\ae}k \cite{groenbaek} characterised commutative semigroups
$G$ for which $\ell^1(G)$ is amenable.
A recent extensive treatise on $\ell^1$-algebras of semigroups
has been written by Dales, Lau and Strauss \cite{dalesls}, which
includes a charaterisation of all semigroups $G$ for which $\ell^1(G)$ is 
amenable.  Biflatness of $\ell^1(S)$, for a semilattice $S$, has recently
been characterised by Choi \cite{choi}.  

\subsection{Preliminaries}
Let $\fA$ be a Banach algebra.  Let $\fA\otimes^\gam\fA$
denote the projective tensor product.  We let
$m:  \fA\otimes^\gam\fA\to\fA$ denote the multiplication map
and we have left and right module actions of $\fA$
on $\fA\otimes^\gam\fA$ given on elementary tensors by
\[
a\mult(b\otimes c)=(ab)\otimes c\aand
(b\otimes c)\mult a=b\otimes (ca).
\]
A {\it bounded approximate diagonal} (b.a.d.) is a bounded
net $(D_\alp)$ in $\fA\otimes^\gam\fA$ such that
$(m(D_\alp))$ is a bounded approximate identity in $\fA$, i.e.\
\begin{equation}\label{eq:appdiag1}
\lim_\alp am(D_\alp)=a\aand\lim_\alp m(D_\alp)a=a\text{ for each }
a\iin\fA
\end{equation}
and $(D_\alp)$ is asymptotically central for the $\fA$-actions, i.e.\
\begin{equation}\label{eq:appdiag2}
\lim_\alp(a\mult D_\alp-D_\alp\mult a)=0\text{ for each }a\iin\fA.
\end{equation}
Following Johnson \cite{johnson}, we will say that a Banach algebra
$\fA$ is {\it amenable} if it admits a b.a.d.
A quantitative feature of amenability was introduced by Johnson
in \cite{johnson95}, for applications to Fourier algebras of finite groups.
The {\it amenability constant} of an amenable Banach algebra $\fA$
is given by
\[
\mathrm{AM}(\fA)=\inf\left\{\sup_\alp\norm{D_\alp}_\gam:(D_\alp)
\text{ is a b.a.d.\ for }\fA\right\}.
\]
The problem of understanding amenable semigroup algebras
in terms of their amenability constants has attracted some attention
\cite{stokke,dalesls}.

We call $\fA$ {\it contractible} if it admits a {\it diagonal}, i.e.\
an element $D\iin\fA\otimes^\gam\fA$ for which
\begin{align}
am(D)=&a=m(D)a\;\aand\label{eq:diag1} \\
a\mult D&=D\mult a\label{eq:diag2}
\end{align}
for each $a\iin\fA$.   Note, in particular, then $\fA$ must be unital
and the norm of the unit is bounded above by $\mathrm{AM}(\fA)$.

If $\fA$ is a finite dimensional amenable Banach algebra, then
$\fA\otimes^\gam\fA$ is a finite dimensional Banach space, so
any b.a.d.\ admits a cluster point $D$.
Since any subnet of a b.a.d.\ is also a b.a.d., the cluster point
must be a diagonal, whence $\fA$ is contractible.

We record the following simple observation.

\begin{unique}\label{prop:unique}
If $\fA$ is a contractible commutative Banach algebra,
then the diagonal is unique.
\end{unique}

\proof We note that $\fA\otimes^\gam\fA$ is a Banach algebra
in an obvious way:  $(a\otimes b)(c\otimes d)=(ac)\otimes(bd)$.
If $D$ is a diagonal, then $(a\otimes b)D=a\mult D\mult b=(ab)\mult D$ for
$a,b\iin\fA$, by commutativity.  Hence if $D'$ is another diagonal
\[
D'D=m(D')\mult D=1\mult D=D
\]
and, similarly, $D'D=DD'=D'$.\endpf

It will also be useful to observe the following.

\begin{tpdiag}\label{prop:tpdiag}
Let $\fA$ and $\fB$ be contractible Banach algebras, 
with respective diagonals $D_\fA$ and $D_\fB$, then
$\fA\otimes^\gam\fB$ has diagonal
\[
D_\fA\otimes D_\fB\in(\fA\otimes^\gam\fA)\otimes^\gam(\fB\otimes^\gam\fB)
\cong(\fA\otimes^\gam\fB)\otimes^\gam(\fA\otimes^\gam\fB).
\]
\end{tpdiag}

\proof  It is simple to check the diagonal axioms (\ref{eq:diag1}) and
(\ref{eq:diag2}). \endpf

\section{Amenability constants for semilattice algebras}
\label{sec:semilattice}

A {\it semilattice} is a commutative semigroup $S$
in which each element is idempotent, i.e.\ if $s\in S$ then $ss=s$.
If $s,t\in S$ we write 
\begin{equation}\label{eq:po}
s\leq t\quad\iff\quad st=s.
\end{equation}
It is clear that this defines a partial order on $S$.  We note
that if $S$ is a finite semilattice, then $o=\prod_{s\in S}s$ is
a minimal element for $S$ with respect to this partial order.
We note that if $S$ has a minimal element, then it is unique.
Also if $S$ has a unit $1$, then $1$ is the maximal element in
$S$.  

A basic example of a semilattice is $\fP(T)$, the set of all
subsets of a set $T$, where we define $\sig\tau=\sig\cap \tau$
for $\sig,\tau\iin\fP(T)$.  The minimal element is $\varnothing$, and
the maximal element is $T$.
We call any subsemilattice of a semilattice $\fP(T)$ a {\it subset
semilattice}.  
This type of semilattice is universal as we have a
semilattice ``Cayley Theorem":  for any semilattice $S$,
the map $s\mapsto\{t\in S:t\leq s\}:S\to\fP(S)$ 
(or $s\mapsto\{t\in S\setdif\{o\}:t\leq s\}:S\to\fP(S\setdif\{o\})$)
is an injective semilattice homomorphism (by which 
$o\mapsto\varnothing$).

For any semilattice $S$ we define
\[
\ell^1(S)=\left\{x=\sum_{s\in S}x(s)\del_s:\text{each }x(s)\in\Cee\aand
\norm{x}_1=\sum_{s\in S}|x(s)|<\infty\right\}
\]
where each $\del_s$ is the usual ``point mass'' function.  Then
$\ell^1(S)$ is a commutative 
Banach algebra under the norm $\norm{\cdot}_1$ with the product 
\[
\left(\sum_{s\in S}x(s)\del_s\right)\con\left(\sum_{t\in S}x(t)\del_t\right)
=\sum_{r\in S}\left(\sum_{st=r}x(s)y(t)\right)\del_r.
\]
In particular we have $\del_s\con\del_t=\del_{st}$.  
We shall consider the Banach space $\ell^\infty(S)$, of bounded functions from 
$S$ to $\Cee$ with supremum norm,
to be an algebra under usual pointwise operations.  The Cayley map,
indicated above, extends to an algebra homomorphism 
$\Sig:\ell^1(S)\to\ell^\infty(S)$, given on each $\del_s$ by
\begin{equation}\label{eq:schutz}
\Sig(\del_s)=\chi_{\{t\in S:t\leq s\}}
\end{equation}
and extended linearly and continuously to all of $\ell^1(S)$.
Here, $\chi_T$ is the indicator function of $T\subset S$.
The map $\Sig$ is called the {\it Sch\"{u}tzenburger map}; 
see \cite[\S 4]{choi}
and references therein.  

We note that if $S$ is finite, then
$\Sig$ is a bijection.  In this case a formula for its inverse is given by
\begin{equation}\label{eq:mobius}
\Sig^{-1}(\chi_s)=\sum_{t\leq s}\mu(t,s)\del_t
\end{equation}
where $\chi_s=\chi_{\{s\}}$ and
$\mu:\{(t,s):S\cross S:t\leq s\}\to\Ree$ is the {\it M\"{o}bius
function} of the partially ordered set $(S,\leq)$ as defined in 
\cite[\S 3.7]{stanley}.  Our computations
in this section will be equivalent to explicitly computating $\mu$, 
though we will never need to know $\mu$ directly.

It follows from \cite[Theorem 10]{duncann}
that $\ell^1(S)$ is amenable if and only if $S$ is finite.
Thus it follows (\ref{eq:diag1}) that $\ell^1(S)$ is unital if $S$ is finite.
If $S$ is unital, then $\del_1$ is the unit
for $\ell^1(S)$.  If $S$ is not unital, the unit is more complicated.
We let $M(S)$ denote the set of maximal elements in $S$ with respect to
the partial ordering (\ref{eq:po}).

\begin{unit}\label{prop:unit}
If $S$ is a finite semilattice then the unit is given by
$u=\sum_{p\in S}u(p)\del_p$ where
\begin{equation}\label{eq:unit}
u(p)=1-\sum_{t>p}u(t)
\end{equation}
for each $p\iin S$ and we adopt the convention that an empty sum is $0$.
Moreover
\begin{equation}\label{eq:unit1}
\sum_{s\in S}u(s)=1.
\end{equation}
\end{unit}

\proof While we have already established existence of the unit above,
let us note that we can gain a very elementary proof of its existence.
Indeed since $\Sig:\ell^1(S)\to\ell^\infty(S)$ is a bijection, 
$u=\Sig^{-1}(\chi_S)$ is the unit for $\ell_1(S)$.

If $p\in S$ then
\[
\del_p=\del_p\con u=\left(\sum_{s\geq p}u(s)\right)\del_p
+\sum_{s<p}\left(\sum_{\substack{t\in S \\ tp=s}}u(t)\right)\del_s.
\]
and thus, inspecting the coefficient of $\del_p$, we obtain (\ref{eq:unit}).
Note that if $p\in M(S)$ the formula above gives $u(p)=1$, and for any 
$s\iin S\setdif M(S)$ we have $\sum_{t\in S,tp=s}u(t)=0$.
Thus, if we select $p\iin M(S)$ we have
\[
\sum_{s\in S}u(s)=u(p)
+\sum_{s<p}\left(\sum_{\substack{t\in S \\ tp=s}}u(t)\right)=1
\]
and thus obtain (\ref{eq:unit1}).  \endpf


We note that if $S$ is a finite 
semilattice then $S\setdif M(S)$ is a subsemilattice, in fact an
ideal, of $S$.   We also note that $S\cross S$
is also a semilattice and the partial order there satisfies
\[
(s,t)\leq(p,q)\quad\iff\quad s\leq p\aand t\leq q.
\]
The following gives an algorithm for computing the diagonal for $\ell^1(S)$.

\begin{algorithm}\label{lem:algorithm}
Let $S$ be a finite semilattice.  Then the diagonal
\[
D=\sum_{(s,t)\in S\times S}d(s,t)\del_s\otimes\del_t
\]
satisfies, for all $(p,q)\iin S\cross S$,

{\bf (a)} $\displaystyle d(p,p)=u(p)
-\sum_{\substack{(s,t)>(p,p) \\ st=p}}d(s,t)$;

{\bf (b)} if $q\not\geq p$, then $\displaystyle d(p,q)=-\sum_{t>q}d(p,t)$
and $\displaystyle d(q,p)=-\sum_{s>q}d(s,p)$; and

{\bf (c)} $d(p,q)=d(q,p)$.

\noindent Thus, each $d(p,q)$ is an integer, and for distinct elements
$p,q\iin M(S)$ we have $d(p,p)=1$ and $d(p,q)=0$.
\end{algorithm}

\proof The equation (\ref{eq:diag1}) gives us
\begin{equation}\label{eq:diag11}
\sum_{p\in S}u(p)\del_p=u=\sum_{(s,t)\in S\times S}d(s,t)\del_{st}
=\sum_{p\in S}\left(\sum_{\substack{(s,t)\in S\times S \\ st=p}}d(s,t)\right)
\del_p
\end{equation}
Since $st=p$ necessitates $(s,t)\geq(p,p)$, we examine the
coefficient of $\del_p$ to find
\begin{equation}\label{eq:plevel}
u(p)
=\sum_{\substack{(s,t)\geq(p,p) \\ st=p}}d(s,t)
\end{equation}
from which we obtain (a).  In particular, if $p\in M(S)$
we obtain an empty sum in (a) and find $d(p,p)=1$.
The equation (\ref{eq:diag2})
implies that $\del_q\mult D=D\mult\del_q$ and hence we obtain
\begin{equation}\label{eq:diag22}
\sum_{(s,t)\in S\times S}d(s,t)\del_{qs}\otimes\del_t
=\sum_{(s,t)\in S\times S}d(s,t)\del_s\otimes\del_{tq}.
\end{equation}
If $q\not\geq p$ then there is no $s\iin S$ for which $qs=p$.
Hence examining the coefficient of $\del_p\otimes\del_q$ and
$\del_q\otimes\del_p$, respectively, in (\ref{eq:diag22}), yields
\begin{equation}\label{eq:dptsum}
0=\sum_{t\geq q}d(p,t)\quad\aand\quad
\sum_{s\geq q}d(s,p)=0.
\end{equation}
Hence we have established (b).  In particular, if $q,p\in M(S)$
we have an empty sum in (b), so $d(p,q)=0$.

We can see for any pair $(p,q)$ with $p\not=q$,
so $p\not\leq q$ or $q\not\leq p$, that $d(p,q)$ is determined
by coefficients $(s,t)>(p,q)$.  
Hence by induction, using the coeficients $d(p,p)$ and $d(p,q)$ for distict
maximal $p,q$ as a base, we obtain (c).
For example, if $q\in M(S\setdif M(S))$, then (b) implies for every
$p>q$ that
\[
d(p,q)=-\sum_{t>q}d(p,t)=-d(p,p)=-1
\]
and, similarly, $d(q,p)=-1$.  

It is clear, form the above induction, that each $d(p,q)$ is an integer.  
\endpf

Let us see how Lemma \ref{lem:algorithm} allows us to compute the diagonal
$D$ of $\ell^1(S)$ for a finite semilattice $S$.

\medskip
{\bf Step 1.} We inductively define
\begin{equation}\label{eq:ideals}
S_0=S,S_1=S\setdif M(S),\dots,S_{k+1}=S_k\setdif M(S_k)
\end{equation}
and we let $n(S)=\min\{k:S_{k+1}=\varnothing\}$, so $S_{n(S)}=\{o\}$
and $S_{n(S)+1}=\varnothing$.

{\bf Step 2.} We label 
$S=\{s_0,s_1,\dots,s_{|S|-1}\}$ in any manner for which 
\[
i\geq j\aand s_i\in S_k\quad\implies\quad s_j\in S_k.
\]
Thus, the elements of $M(S_k)$ comprise the last part
of the list of $S_k$ for $k=1,\dots,n(S)$.  In particular,
$s_0=o$ and $s_{|S|-1}\in M(S)$.

{\bf Step 3.} The diagonal $D$ will be represented by an
$|S|\cross|S|$ matrix $[D]=[d(s_i,s_j)]$.  The lower rightmost
corner will be the $|M(S)|\cross|M(S)|$ identity matrix.
We can then proceed, using formulas (b) and
(a) from the lemma above, to compute the remaining
entries of the lower rightmost
$(|M(S)|+1)\cross(|M(S)|+1)$ corner of $[D]$, etc., until we are done.

\medskip
In order to describe certain semilattices $S$,
we define the {\it semilattice graph} 
$\Gamma(S)=(S,e(S))$, where the vertex set is $S$ and the edge set
is given by ordered pairs
\[
e(S)=\{(s,t)\in S\cross S:s>t\text{ and there is no }r\iin S\text{ for which } 
s>r>t\}.
\]
To picture such a graph for a finite semilattice $S$
it is helpful to describe levels.  
Let $S_0,S_1,\dots,S_{n(S)}$ be the sequence of ideals of $S$ given in
(\ref{eq:ideals}).  For $s\iin S$ we let the
{\it level} of $s$ be given by
\[
\lam(s)=n(S)-k\wwhere s\in M(S_k). 
\]
Note that for the power set semilattice $\fP(T)$, $\lam(\sig)=|\sig|$,
the cardinality of $\sig$.  However, this relation need not hold for 
a subsemilattice of $\fP(T)$, as is evident from the Example
\ref{ex:flatsem}, below.  A 6-element, 4-level semilattice is
illustrated in (\ref{ex:unrsl}).

We apply this algorithm to obtain the following examples.
We denote, for a finite semilattice $S$, the amenability constant
\[
\mathrm{AM}(S)=\mathrm{AM}(\ell^1(S))=\norm{D}_1=
\sum_{(s,t)\in S\times S}|d(s,t)|
\]
where we recall the well-known isometric identification
$\ell^1(S)\otimes^\gam\ell^1(S)\cong\ell^1(S\cross S)$.

\begin{linear}\label{ex:linear}
Let $L_n=\{0,1,2,\dots,n\}$ be a ``linear" semilattice with
operation $st=s\wedge t=\min\{s,t\}$.  Then we obtain diagonal
with $(n+1)\cross(n+1)$ matrix
\[
[D]=\begin{bmatrix}
\phantom{-}2 & -1 & \hdots & \phantom{-}0 & \phantom{-}0 \\
-1 & \phantom{-}2 & \ddots & \phantom{-}0 & \phantom{-}0 \\
\phantom{-}\vdots & \ddots & \ddots & \ddots & \phantom{-}\vdots \\
\phantom{-}0 & \phantom{-}0 & \ddots & \phantom{-}2 & -1 \\
\phantom{-}0 & \phantom{-}0 & \hdots &  -1 & \phantom{-}1
\end{bmatrix}.
\]
Hence $\mathrm{AM}(L_n)=4n+1$.
\end{linear}

\begin{flatsem}\label{ex:flatsem}
Let $F_n=\{o,s_1,\dots,s_n\}$ be the $n+1$ element ``flat'' semilattice with
multiplications $s_is_j=o$ if $i\not=j$.  Then we obtain unit
\[
u=\del_{s_1}+\dots+\del_{s_n}+(1-n)\del_o
\]
and diagonal with $(n+1)\cross(n+1)$-matrix
\[
[D]=\begin{bmatrix}
n+1 & -1 & -1 & \hdots & -1 \\
-1 & \phantom{-}1 & \phantom{-}0 & \hdots & \phantom{-}0 \\
-1 & \phantom{-}0 & \phantom{-}1 & \ddots & \phantom{-}\vdots \\
\phantom{-}\vdots & \phantom{-}\vdots & \ddots & \ddots & \phantom{-}0 \\
-1 & \phantom{-}0 & \hdots & \phantom{-}0 & \phantom{-}1
\end{bmatrix}.
\]
Hence $\mathrm{AM}(F_n)=4n+1$.
\end{flatsem}

\begin{uflatsem}\label{ex:uflatsem}
Let $F_n^1=\{o,s_1,\dots,s_n,1\}$ be 
the unitasation of $F_n$, above.
Then we obtain diagonal with 
$(n+2)\cross(n+2)$ matrix
\[
[D]=\begin{bmatrix}
n^2-n+2 & -n & \hdots & -n & n-1 \\
-n & \phantom{-}2 & \hdots & \phantom{-}1 & -1 \\ 
\vdots & \phantom{-}\vdots & \ddots & \phantom{-}\vdots & \vdots \\
-n & \phantom{-}1 & \hdots & \phantom{-}2 & -1 \\
n-1 & -1 & \hdots & -1 & \phantom{-}1\end{bmatrix}.
\]
Hence $\mathrm{AM}(F_n^1)=4n^2+4n+1$.
\end{uflatsem}

The next example is less direct than the previous ones, so we offer a proof.

\begin{powset}\label{ex:powset}
Let $P_n=\fP(\{1,\dots,n\})$ with multiplication $st=s\cap t$.
Then the diagonal $D$ has $2^n\cross 2^n$ matrix which is, up
to permutative similarity, the Kronecker product
\[
\begin{bmatrix} \phantom{-}2 & -1 \\ -1 & \phantom{-}1\end{bmatrix}\otimes\dots
\otimes\begin{bmatrix} \phantom{-}2 & -1 \\ -1 & \phantom{-}1
\end{bmatrix}\;(n\text{ times}). 
\]
Hence $\mathrm{AM}(P_n)=5^n$.
\end{powset}

\proof If $s\in P_n$ let $\chi_s:\{1,\dots,n\}\to\{0,1\}=L_1$
be its indicator function.  It is easily verified that the map 
$s\mapsto\chi_s:P_n\to L_1^n$ is a semilattice isomorphism.
Thus there is an isometric identification $\ell^1(P_n)\cong
\ell^1(L_1)\otimes^\gam\dots\otimes^\gam\ell^1(L_1)$.
Then it follows from Proposition \ref{prop:tpdiag} above that 
$D=D_1\otimes\dots\otimes D_1$
where $D_1$ is the diagonal for $\ell^1(L_1)$, which, by the algorithm
has matrix
\[
[D_1]=\begin{bmatrix} \phantom{-}2 & -1 \\ -1 & \phantom{-}1\end{bmatrix}.
\]
The amenability constant $\mathrm{AM}(P_n)$ can be easily computed
by induction.  \endpf

We have the following summary result.

\begin{amenconst}\label{theo:amenconst}
If $S$ is a finite semilattice, then $\mathrm{AM}(S)=4n+1$
for some integer $n\geq 0$.  All such numbers are achieved.
\end{amenconst}

\proof We first establish that for $p\iin S$, $d(p,p)\geq 0$.
This does not seem obvious from Lemma \ref{lem:algorithm}.  We use
a calculation from \cite[\S 3]{choi} which exploits the M\"{o}bius function.
We have that $\Sig:\ell^1(S)\to\ell^\infty(S)$ is invertible
and $\til{D}=\sum_{r\in S}\chi_r\otimes\chi_r$ is the diagonal for
$\ell^\infty(S)$.  Thus, using (\ref{eq:mobius}), we have that 
\begin{align*}
D=\Sig^{-1}\otimes\Sig^{-1}(\til{D})
&=\sum_{r\in S}\left(\sum_{s\in S}\til{\mu}(s,r)\del_s\right)
\otimes\left(\sum_{t\in S}\til{\mu}(t,r)\del_t\right) \\
&=\sum_{(s,t)\in S\times S}\left(\sum_{r\in S}\til{\mu}(s,r)
\til{\mu}(t,r)\right)\del_s\otimes\del_t
\end{align*}
is the diagonal for $\ell^1(S)$,
where $\til{\mu}(s,t)=\mu(s,t)$ if $s\leq t$ and $\til{\mu}(s,t)=0$, otherwise.
Inspecting the coeficient of $\del_p\otimes\del_p$ we obtain
\begin{equation}\label{eq:dpppos}
d(p,p)=\sum_{r\in S}\til{\mu}(p,r)^2\geq 1> 0
\end{equation}
since $\til{\mu}(p,p)=\mu(p,p)=1$ by \cite[\S 3.7]{stanley}.
We now observe, using (\ref{eq:plevel}) and then (\ref{eq:unit1}), 
that
\[
\sum_{(s,t)\in S\times S}d(s,t)=
\sum_{p\in S}\sum_{\substack{(s,t)\in S\times S \\ st=p}}d(s,t)
=\sum_{p\in S}u(p)=1.
\]
By symmetry, if $p\not=q$ then $|d(p,q)|+|d(q,p)|\equiv d(p,q)+d(q,p)\mod 4$.
Hence we have
\[
\mathrm{AM}(S)\equiv\sum_{(s,t)\in S\times S}|d(s,t)|
\equiv\sum_{(s,t)\in S\times S}d(s,t)\equiv 1\mod 4.
\]
Finally, Examples \ref{ex:linear} and \ref{ex:flatsem} provide us with
semilattices admitting amenability constants $4n+1$, for each integer 
$n\geq 0$.  \endpf

We now gain a crude lower bound for $\mathrm{AM}(S)$ which we will
require for Proposition \ref{prop:dalesq}.

\begin{amenconst1}\label{cor:amenconst1}
For any finite semilattice $S$ we have $\mathrm{AM}(S)\geq 2|S|-1$.
\end{amenconst1}

\proof  We have from (\ref{eq:dpppos}) that $d(p,p)\geq 1$
for each $p\iin S$.  It then follows from (\ref{eq:dptsum})
that for $p>o$ we have $\sum_{t\geq o}d(p,t)=0$ from which we obtain
$\sum_{t\not=p}|d(p,t)|\geq 1$.  It then follows that
\begin{align*}
\mathrm{AM}(S)=\sum_{(s,t)\in S\times S}|d(s,t)|
&\geq d(o,o)+\sum_{p>0}\left(d(p,p)+\sum_{t\not=p}|d(p,t)|\right) \\
&\geq 1+(|S|-1)2
\end{align*}
and we are done.  \endpf

We note that if $S$ is unital, then for $p<1$, $u(p)=0$ and since
$d(s,t)=d(t,s)$ for $(s,t)>(p,p)$ we find from Lemma \ref{lem:algorithm}
(a) that $d(p,p)$ is even; in particular $d(p,p)\geq 2$.  The proof above
may be adapted to show $\mathrm{AM}(S)\geq 4|S|-3$, in this case.
We conjecture the estimate $\mathrm{AM}(S)\geq 4|S|-3$ holds for any
finite semilattice $S$.

\section{Banach algebras graded over semilattices}
\label{sec:graded}

A Banach algebra $\fA$ is {\it graded} over a semigroup 
$S$ if we have closed subspaces $\fA_s$ for each $s\iin S$ such that
\[
\fA=\ell^1\text{-}\bigoplus_{s\in S}\fA_s\aand \fA_s\fA_t\subset\fA_{st}
\ffor s,t\iin S.
\]
We will be interested strictly in the case where $S$ is a finite 
semilattice.  Notice in this case each $\fA_s$ is a closed subalgebra of $\fA$.
The next proposition can be proved by a straightforward adaptation
of the proof of \cite[Proposition 3.1]{rundes2}.  However, we offer another proof.

\begin{gradedalga}\label{prop:gradedalga}
Let $S$ be a finite semilattice and $\fA$ be graded over $S$.
Then $\fA$ is amenable if and only if each $\fA_s$ is amenable.
\end{gradedalga}

\proof Suppose $\fA$ is amenable.
If $s\in S$, then $\fA^s=\bigoplus_{t\leq s}\fA_t$
is an ideal in $\fA$ which is complemented and hence an amenable
Banach algebra (see \cite[Theorem 2.3.7]{rundeB}, for example).  
It is easy the check that
the projection $\pi_s:\fA^s\to\fA_s$ is a quotient homomorphism.
Hence it follows that if $(D^s_\alp)$ is an approximate diagonal for
$\fA^s$ then $\bigl(\pi_s\otimes\pi_s(D^s_\alp)\bigr)$ is an approximate 
diagonal for  $\fA_s$.  (This is quotient argument is noted in 
\cite[Corollary 2.3.2]{rundeB} and \cite[Proposition 2.5]{dalesls}.)

Now suppose that each $\fA_s$ is amenable.
Let $S_0,S_1,\dots, S_{n(S)}$ be the sequence of ideals from
(\ref{eq:ideals}).
For each $n=0,1,\dots,n(S)$ we set $\fA_n=\bigoplus_{s\in S_n}\fA_s$
and observe, for each $n=0,1,\dots,n(S)-1$, that we have an isometrically
isomorphic identification
\[
\fA_n/\fA_{n+1}=\ell^1\text{-}\!\!\!\bigoplus_{s\in M(S_n)}\fA_s
\]
where multiplication in the latter is pointwise, i.e.\ $\fA_s\fA_t=\{0\}$
if $s\not=t\iin M(S_n)$.  The pointwise algebra 
$\ell^1\text{-}\bigoplus_{s\in M(S_n)}\fA_s$ is amenable as each $\fA_s$ is
amenable; if $(D_{s,\alp})$ is a bounded approximate diagonal for
each $\fA_s$, then in
\[
\left(\ell^1\text{-}\!\!\!\bigoplus_{s\in M(S_n)}\fA_s\right)\otimes^\gam
\left(\ell^1\text{-}\!\!\!\bigoplus_{s\in M(S_n)}\fA_s\right)
\cong\;\;\;\ell^1\text{-}\!\!\!\!\!\!\!\!\!\!\!\!\!\!
\bigoplus_{(s,t)\in M(S_n)\times M(S_n)}\fA_s\otimes^\gam\fA_t
\]
the net of elements $D_\alp=\sum_{s\in M(S_n)}D_{s,\alp}$ is
an approximate diagonal.
Thus if  $\fA_{n+1}$ is amenable, then $\fA_n$ must be too
by \cite[Theorem 2.3.10]{rundeB}.  The algebra $\fA_{n(S)}=\fA_o$ is amenable, 
and hence we may finish by an obvious induction.  \endpf

In the computations which follow, we will require one of the
following {\it linking assumptions} which are very natural for our examples.

\smallskip
\parbox{4.5in}{{\bf (LA1)}  
For each $s\iin S$ there is a bounded approximate identity
$(u_{s,\alp})_\alp$ in $\fA_s$, such that for each $t\leq s$ and
$a_t\in\fA_t$ we have $\lim_\alp u_{s,\alp} a_t=a_t=\lim_\alp a_tu_{s,\alp}$.}

\smallskip
\parbox{4.5in}{\noindent {\bf (LA2)}  
For each $s\in S$ there is a contractive character 
$\chi_s:\fA_s\to\Cee$ such that for each $s,t\iin S$,
$a_s\in\fA_s$ and $a_t\in\fA_t$, we have
$\chi_{st}(a_sb_t)=\chi_s(a_s)\chi_t(a_t)$.}

\smallskip
\noindent Notice that in (LA1), each $(u_{s,\alp})_\alp$ is a 
bounded approximate identity for 
$\fA^s=\ell^1\text{-}\bigoplus_{t\leq s}\fA_t$.  
Thus since
$\fA^s$ is an $\fA_s$-module, 
Cohen's factorisation theorem \cite[32.22]{hewittrII} tells us that
\begin{equation}\label{eq:cohen}
\text{for each }a\iin\fA^s\text{ there is }v_s\in\fA_s\aand a'\iin\fA^s
\text{ such that }a=v_sa'.
\end{equation}
There is a right factorisation analogue, and the result also holds
on each $\fA_s$ module $\fA_t$, where $t\leq s$.
We note that (LA2) is equivalent to having a contractive
character $\chi:\fA\to\Cee$ such that $\chi|_{\fA_s}=\chi_s$ for each $s$.

We note that many natural Banach algebras, graded over semilattices, 
which arise in harmonic analysis, satisfy (LA2).  However, (LA1)
can be used whenever each component algebra $\fA_s$ admits no characters.
For example, if we have a (finite unital) semilattice $S$, a family
of algebras $\{\fA_s\}_{s\in S}$ each having no characters, and a system
$\{\eta^s_t:s,t\in S, s\geq t\}$ of homomorphisms, we can make
$\ell^1\text{-}\bigoplus_{s\in S}\fA_s$ into a Banach algebra by
setting $a_sa_t=\eta^s_{st}(a_s)\eta^t_{st}(a_t)$ for $a_s\iin \fA_s$
and $a_t\iin\fA_t$.  (This construction is analagous to that of the 
Clifford semigroup algebras which will be presented in Section
\ref{ssec:clifford}, below.)  

This brings us to the main result of this article.

\begin{gradedalg}\label{theo:gradedalg}
Let $\fA$ be a Banach algebra graded over a finite semilattice $S$
such that each $\fA_s$ is amenable.
If we have either that (LA1) holds, or (LA2) holds,
then $\mathrm{AM}(\fA)\geq\mathrm{AM}(S)$.
\end{gradedalg}

\proof  $\fA$ is amenable by the proposition above.

Let us suppose (LA1) holds.  We let for each $p\iin S$,
$\pi_p:\fA\to\fA_p$ the contractive projection.  We define
for $a,b\in\fA$, $\pi_p(a\otimes b)=\pi_p(a)\otimes b$ and
$(a\otimes b)\pi_p=a\otimes\pi_p(b)$.  Clearly these
actions extend linearly and continuously to define
$\pi_pD$ and $D\pi_p$ for any $D\in\fA\otimes^\gam\fA$.

We let $(D_\alp)$ be a bounded approximate diagonal for 
$\fA$ and 
\[
D=\sum_{(s,t)\in S\times S}d(s,t)\del_s\otimes\del_t
\]  
be the unique diagonal for $\ell^1(S)$.
We will prove that for $p,q\in S$ and $a\in\fA^p$, $b\in\fA^q$
that
\[
\lim_\alp am(\pi_pD_\alp\pi_q)b=d(p,q)ab.\tag{$\bigstar$}
\]
This requires induction and we will need some preliminary steps.

Suppose that $q\not=p$ in $S$, say $q\not\geq p$.
If $v_q\in\fA_q$ then (\ref{eq:appdiag2}) implies that
\begin{equation}\label{eq:capdiagvanish}
\lim_\alp \pi_p(D_\alp\mult v_q)\pi_q
=\lim_\alp\pi_p(v_q\mult D_\alp)\pi_q=0.
\end{equation}
We note that on an elementary tensor in $\fA\otimes\fA$ we have
\begin{equation}\label{eq:tensform1}
m(\pi_p(a\otimes b\mult v_q)\pi_q)=\sum_{t\geq q}\pi_p(a)\pi_t(b)v_q
=\sum_{t\geq q}m(\pi_p(a\otimes b)\pi_t)v_q
\end{equation}
Now if $b\in\fA^q$ we find $v_q\in\fA_q$
and $b'\iin\fA^q$ such that $b=v_qb'$ by (\ref{eq:cohen}).
We then have, in analogy to Lemma \ref{lem:algorithm} (b), using 
(\ref{eq:capdiagvanish}) and (\ref{eq:tensform1})
\[
\lim_\alp\sum_{t\geq q}m(\pi_pD_\alp\pi_t)b
=\lim_\alp m\bigl(\pi_p(D_\alp\mult v_q)\pi_q\bigr)b'=0\mult b'=0.
\tag{b$_1$'}
\]
Similarly we see
\[
\lim_\alp\sum_{s\geq q}bm(\pi_sD_\alp\pi_p)=0.
\tag{b$_2$'}
\]
Note that if $p,q\in M(S)$ with $p\not=q$, then then (b$_1$')
takes the form
\[
\lim_\alp m(\pi_pD_\alp\pi_q)b=0=d(p,q)b
\]
and a simlar version holds for (b$_2$').  Thus ($\bigstar$)
holds in this case.

Now we show that for $p\in S$ and $b\iin \fA^p$
that 
\begin{equation}\label{eq:approxuatp}
\lim_\alp\pi_p\bigl(m(D_\alp)\bigr)b=u(p)b
\end{equation}
where $u=\sum_{p\in S}u(p)\del_p$ is the unit for $\ell^1(S)$.
By (\ref{eq:cohen}) there are $v_p\iin\fA_p$ and $b'\iin\fA^p$
such that $b=v_pb'$.  We have that
\begin{align*}
v_p&=\lim_\alp m(D_\alp)v_p
=\lim_\alp\sum_{s\in S}\pi_s\bigl(m(D_\alp)\bigr)v_p \\
&=\lim_\alp\left[\sum_{s\geq p}\pi_s\bigl(m(D_\alp)\bigr)v_p
+\sum_{s\not\geq p}\pi_s\bigl(m(D_\alp)\bigr)v_p\right]
\end{align*}
from which it follows that
\[
\lim_\alp\sum_{s\geq p}\pi_s\bigl(m(D_\alp)\bigr)v_p
=\lim_\alp\pi_p\bigl(m(D_\alp)v_p\bigr)=v_p
\]
and hence
\begin{equation}\label{eq:approxuatp1}
\lim_\alp\sum_{s\geq p}\pi_s\bigl(m(D_\alp)\bigr)b
=\lim_\alp\sum_{s\geq p}\pi_s\bigl(m(D_\alp)\bigr)v_pb'=v_pb'=b.
\end{equation}
In particular, if $p\in M(S)$, then 
\[
\lim_\alp\pi_p\bigl(m(D_\alp)\bigr)b=b=u(p)b.
\]
Then the equation
(\ref{eq:approxuatp}) follows inductively from (\ref{eq:approxuatp1})
and (\ref{eq:unit}), using the case of maximal $p$ as a base.

Now we establish an analogue of Lemma \ref{lem:algorithm} (a).  
For an elementary tensor $a\otimes b$
in $\fA\otimes\fA$, we have 
\begin{equation}\label{eq:tensorform}
\pi_p(ab)=\sum_{\substack{(s,t)\geq (p,p) \\ st=p}}\pi_s(a)\pi_t(b)
=\sum_{\substack{(s,t)\geq (p,p) \\ st=p}}m\bigl(\pi_s(a\otimes b)\pi_t\bigr).
\end{equation}
It then follows from (\ref{eq:approxuatp}) and (\ref{eq:tensorform})
that for $b\in\fA^p$
\[
u(p)b=\lim_\alp\sum_{\substack{(s,t)\geq (p,p) \\ st=p}}
m(\pi_sD_\alp\pi_t)b. \tag{a'}
\]
Note that if $p\in M(S)$, then by Proposition \ref{prop:unit}, (a')
becomes 
\[
d(p,p)b=b=\lim_\alp m(\pi_pD_\alp\pi_p)b.
\]
Thus ($\bigstar$) holds in this case.

We now prove ($\bigstar$) by induction on pairs $(p,q)\iin S\cross S$ 
with pairs $(p,q)\in M(S)\cross M(S)$ as a base.  
If $p\in S$, the induction hypothesis is that for $a,b\in\fA^p$
\[
\lim_\alp am(\pi_sD_\alp\pi_t)b=d(s,t)ab\ffor (s,t)>(p,p)\wwith st=p.
\]
Notice that in the hypothesis above we have $\fA^p\subset\fA^s\cap\fA^t$,
and, moreover, either $t\not\geq s$ or $s\not\geq t$.
But then it follows from (a') 
and Lemma \ref{lem:algorithm} (a) that
\[
\lim_\alp am(\pi_pD_\alp\pi_p)b=\left[u(p)-\sum_{\substack{(s,t)>(p,p) \\
st=p}}d(s,t)\right]ab=d(p,p)ab
\]
which establishes ($\bigstar$) in this case. 
Also, if $q\not=p$, say $q\not\geq p$, then for $a\iin\fA^p$ and
$b\iin\fA^q$ the induction hypothesis is that
\[
\lim_\alp am(\pi_pD_\alp\pi_t)b=d(p,t)ab\ffor t>q.
\]
Combining this with (b$_1$') and Lemma \ref{lem:algorithm} (b)
we obtain the equation ($\bigstar$) for this case.  
We can use (b$_2$') in place of (b$_1$') above, to acheive ($\bigstar$)
with $p$ and $q$ interchanged.

We now use ($\bigstar$) to finish the proof.  Let for $p,q\iin S$
\[
\eta(p,q)=\sup_{a\in\fA^p,b\in\fA^q}\frac{\norm{ab}}{\norm{a}\norm{b}}.
\]
We note that our assumption (LA1) provides that $\eta(p,q)>0$.
For $\eps>0$ let $a_\eps\iin\fA^p$ and $b_\eps\iin\fA^q$ be
so $\frac{\norm{a_\eps b_\eps}}{\norm{a_\eps}\norm{b_\eps}}\geq
(1-\eps)\eta(p,q)$.  Then by ($\bigstar$) we have
\begin{align*}
|d(p,q)|\norm{a_\eps b_\eps}
&=\lim_\alp\norm{a_\eps m(\pi_pD_\alp\pi_q)b_\eps}
\leq\liminf_\alp\norm{a_\eps m(\pi_pD_\alp\pi_q)}\norm{b_\eps} \\
&\leq\liminf_\alp\frac{\norm{a_\eps m(\pi_pD_\alp\pi_q)}}
{\norm{a_\eps}\norm{m(\pi_pD_\alp\pi_q)}}
\norm{a_\eps}\norm{b_\eps}\norm{m(\pi_pD_\alp\pi_q)} \\
&\leq \eta(p,q) \norm{a_\eps}\norm{b_\eps}
\liminf_\alp\norm{m(\pi_pD_\alp\pi_q)}
\end{align*}
which implies
\[
(1-\eps)|d(p,q)|\leq\liminf_\alp\norm{m(\pi_pD_\alp\pi_q)}
\leq\liminf_\alp\norm{\pi_pD_\alp\pi_q}_\gam.
\]
Thus 
\begin{align*}
\mathrm{AM}(S)&=\sum_{(p,q)\in S\times S}|d(p,q)|
\leq\sum_{(p,q)\in S\times S}\liminf_\alp\norm{\pi_pD_\alp\pi_q}_\gam \\
&\leq \liminf_\alp\sum_{(p,q)\in S\times S}\norm{\pi_pD_\alp\pi_q}_\gam 
\overset{(\dagger)}=\liminf_\alp \norm{D_\alp}\leq\sup_\alp\norm{D_\alp}_\gam
\end{align*}
where the equality $(\dagger)$ holds because of the isometric identification
\[
\fA\otimes^\gam\fA=\left(
\ell^1\text{-}\bigoplus_{s\in S}\fA_s\right)\otimes^\gam
\left(\ell^1\text{-}\bigoplus_{t\in S}\fA_t\right)
\cong\ell^1\text{-}\!\!\!\!\!\bigoplus_{(s,t)\in S\times S}
\fA_s\otimes^\gam\fA_t.
\]
Thus we have finished the case where we assumed (LA1).

Now suppose we have (LA2).  The map
\[
\Pi:\fA\to\ell^1(S),\quad \Pi(a)=\sum_{s\in S}\chi_s\bigl(\pi_s(a)\bigr)\del_s
\]
is a contractive homomorphism.  Hence it follows that
if $(D_\alp)$ is a bounded approximate diagonal for $\fA$
then $\bigl(\Pi(D_\alp)\bigr)$ is an approximate diagonal for
$\ell^1(S)$.   Thus the limit point, i.e.\ unique cluster point, $D$ of
$\bigl(\Pi(D_\alp)\bigr)$ satisfies $\norm{D}_\gam=\mathrm{AM}(S)$,
whence $\sup_\alp\norm{D_\alp}_\gam\geq\lim_\alp\norm{\Pi(D_\alp)}_\gam
\geq\mathrm{AM}(S)$.  \endpf

It might seem plausible that in the situation of the theorem above, if
it were the case that $\mathrm{AM}(\fA_s)=1$, for each $s$, then
$\mathrm{AM}(\fA)=\mathrm{AM}(S)$.  Indeed this phenomenon was observed
for $S=L_1$, in a special case in \cite[Theorem 2.3]{rundes2}.
However this does not seem to hold in general, as we shall see below.

\subsection{Clifford semigroup algebras}\label{ssec:clifford}
Let $S$ be a semilattice, and for each $s\iin S$ suppose we have
a group $G_s$, and for each $t\leq s$ a homomorphism $\eta^s_t:G_s\to G_t$
such that for $r\geq s\geq t\iin S$ we have
\[
\eta^s_s=\mathrm{id}_{G_s}\quad\aand\quad
\eta^r_s\comp\eta^s_t=\eta^r_t
\]
then $G=\bigsqcup_{s\in S}G_s$ (disjoint union) admits a semigroup
operation given by
\[
x_sy_t=\eta^s_{st}(x_s)\eta^t_{st}(y_t)
\]
for $x_s\iin G_s$ and $y_t\iin G_t$.  It is straightforward to check
that $G$ is a semigroup, and is called a {\it Clifford semigroup}, as such
a semigroup was first described in \cite{clifford}.  We note that
the set of idempotents $E(G)$ is $\{e_s\}_{s\in S}$, where $e_s$ is the
neutral element of $G_s$, and $E(G)$ is a subsemigroup, isomorphic to $S$.
It is clear that
\[
\ell^1(G)=\ell^1\text{-}\bigoplus_{s\in S}\ell^1(G_s)
\]
and that $\ell^1(G)$ is thus graded over $S$.  Note that
$\ell^1(G)$ satisfies (LA1) by design,
and satisfies (LA2) where the augmentation
character is used on each $\ell^1(G_s)$.
As with semilattices we will write $\mathrm{AM}(G)=\mathrm{AM}(\ell^1(G))$

Consider the semilattice $S=\{o,s_1,s_2,s_3,s_4,1\}$ whose
graph is given below.
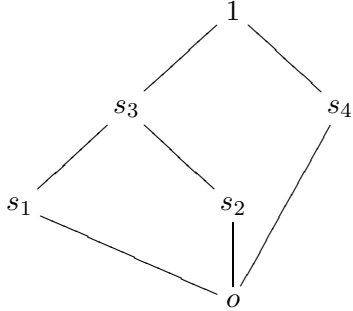
\begin{equation}\label{ex:unrsl}
\xymatrix{
& &  1\ar@{-}[dl]\ar@{-}[dr] & \\
& s_3\ar@{-}[dl]\ar@{-}[dr] &  & s_4\ar@{-}[ddl] \\
s_1\ar@{-}[drr] & & s_2\ar@{-}[d] & & \\
& & o & & 
}
\end{equation}
Using the algorithm following
Lemma \ref{lem:algorithm}, with the semilattice ordered as presented,
we obtain diagonal
$D$ with matrix
\begin{equation}\label{eq:unrslm}
[D]=\begin{bmatrix}
 \phantom{-}6 & -2 & -2 &  \phantom{-}0 & -2 &  \phantom{-}1 \\
-2 &  \phantom{-}2 &  \phantom{-}1 & -1 &  \phantom{-}0 &  \phantom{-}0 \\
-2 &  \phantom{-}1 &  \phantom{-}2 & -1 &  \phantom{-}0 &  \phantom{-}0 \\
 \phantom{-}0 & -1 & -1 &  \phantom{-}2 &  \phantom{-}1 & -1 \\
-2 &  \phantom{-}0 &  \phantom{-}0  &  \phantom{-}1 &  \phantom{-}2 & -1 \\
\phantom{-} 1 &  \phantom{-}0 &  \phantom{-}0  & -1 & -1 &  \phantom{-}1
\end{bmatrix}. 
\end{equation}
Thus we obtain amenability constant $\mathrm{AM}(S)=41$.

Now let $n\geq 2$ be an integer and
$G_n$ be the Clifford semigroup graded over $S$ for which
\[
G_{n,s_3}=\{e_3,a,\dots, a^{n-1}\}\quad\aand\quad G_{n,s_i}=\{e_i\}
\text{ for all }i\not=3
\]
and all connecting homomorphisms are trivial.
Here, $\{e_3,a,\dots,a^{n-1}\}$ is a cyclic group, and each other
$\{e_i\}$ is the trivial group.  This is a finite dimensional
commutative amenable algebra, and hence admits a
unique diagonal by Proposition \ref{prop:unique}.
It is straightforward to verify that if we order the semigroup
$\{o,e_1,e_2,e_3,a,\dots,a^{n-1},e_4,1\}$ we obtain matrix for the diagonal
\[
\begin{bmatrix}
\phantom{-}6  & -2  & -2 &  (1-n)/n & \phantom{-}1/n &
\hdots &  \phantom{-}1/n  & -2 &  \phantom{-}1 \\
-2   &  \phantom{-}2 &  \phantom{-}1 & -1/n &  -1/n &
\hdots & -1/n  &  \phantom{-}0 &  \phantom{-}0 \\
-2 &  \phantom{-}1 &  \phantom{-}2 & -1/n &  -1/n &
\hdots &  -1/n & \phantom{-}0 &  \phantom{-}0 \\
(1-n)/n  & -1/n & -1/n &  (n+1)/n & \phantom{-}0  & \hdots &  
\phantom{-}0 &  \phantom{-}1 & -1 \\
\phantom{-}1/n & -1/n & -1/n &  
\phantom{-}0  & & \iddots & \phantom{-}1/n  
& \phantom{-}0 & \phantom{-}0 \\
\phantom{-}\vdots & \phantom{-}\vdots & \phantom{-}\vdots &
\phantom{-}\vdots & \phantom{-}\iddots & \iddots &  
& \phantom{-}\vdots & \phantom{-}\vdots \\
\phantom{-}1/n & -1/n & -1/n & \phantom{-}0 & 
\phantom{-}1/n &  & \phantom{-}0 & \phantom{-}0 &
\phantom{-}0 \\
-2 &  \phantom{-}0 &  \phantom{-}0  &  \phantom{-}1 & \phantom{-}0  
& \hdots & \phantom{-}0 &  \phantom{-}2 & -1 \\
\phantom{-}1 &  \phantom{-}0 &  \phantom{-}0  & -1 & \phantom{-}0  
& \hdots & \phantom{-}0 & -1 &  \phantom{-}1
\end{bmatrix}.
\]
Notice that values in positions $(o,e_3),\dots,
(o,a^{n-1})$ sum to $0$, the value in the $(o,s_3)$ position 
in (\ref{eq:unrslm}) above.  Similar results holds for all submatirices
with indicies from $\{e_3,a,\dots,a^{n-1}\}$.
Summing absolute values of all entries in the matrix 
we obtain amenability constant $\mathrm{AM}(G_n)=41+4(n-1)/n$.
Thus 
\[
\mathrm{AM}(G_n)=41+4\frac{n-1}{n}>41=\mathrm{AM}(S).
\]
The constant $\mathrm{AM}(G_2)=43$ is the smallest amenability constant we can 
find for an commutative semigroup which is not of the form $4n+1$.




\subsection{Algebras graded over linear semilattices}
We note that if $G$ is a finite Clifford semigroup, graded over a linear 
semilattice $L_n$, then
$\mathrm{AM}(G)=\mathrm{AM}(L_n)=4n+1$.  Indeed, this holds
more generally, by the following proposition.

\begin{gradedlinear}\label{prop:gradedlinear}
If $\fA=\ell_1\text{-}\bigoplus_{k\in L_n}\fA_k$ is a graded
Banach algebra which satisfies (LA1), and
$\fA_k$ is contractible with
$\mathrm{AM}(\fA_k)=1$ for each $k\iin L_n$,
then $\mathrm{AM}(\fA)=4n+1$.
\end{gradedlinear}

\proof  We have from Theorem \ref{theo:gradedalg} that $\mathrm{AM}(\fA)
\geq\mathrm{AM}(L_n)=4n+1$, hence it suffices to exhibit a diagonal
$D$ with $\norm{D}_\gam\leq 4n+1$.  We will show that such $D$ exists
by induction.

Write $L_n=\{0,1,\dots,n\}$.  We identify $L_k$ as an ideal of $L_n$ 
for each $k=0,1,\dots,n-1$ in the usual way.
Let us note that if $(u_{k,\alp})$ is a bounded approximate identity
for $\fA_k$, which satisfies (LA1), then the unit $e_k$ of $\fA_k$
is the limit point of $(u_{k,\alp})$, and hence $e_k$ is the unit
for $\fA^k=\ell^1\text{-}\bigoplus_{j\in L_k}\fA_j$.  Note, moreover,
that the assumption that $\mathrm{AM}(\fA_k)=1$ forces $\norm{e_k}=1$.

Let $\eps>0$.
Suppose for $k<n$ we have a diagonal $D^k$ for $\fA^k$ with 
$\norm{D^k}_\gam<4k+1+\eps$.  For $k=0$, such a diagonal
exists as $\mathrm{AM}(\fA_0)=1$.  We let
\[
D_{k+1}=\sum_{i=1}^\infty a_i\otimes b_i,\quad a_i,b_i\in\fA_{k+1}
\]
be a diagonal for $\fA_{k+1}$ with $\norm{D_{k+1}}_\gam\leq
\sum_{i=1}^\infty\norm{a_i}\norm{b_i}<1+\eps$.
We then set
\[
D^{k+1}
=\sum_{i=1}^\infty a_i\mult\bigl((e_{k+1}-e_k)\otimes(e_{k+1}-e_k)
+D^k\bigr)\mult b_i.
\]
Clearly
\[
\norm{D^{k+1}}_\gam\leq (4+(4k+1+\eps))(1+\eps)=4(k+1)+1+O(\eps).
\]
Applying the multiplication map, and noting that $m(D^k)=e_k$, we have
\begin{align*}
m(D^{k+1})&=\sum_{i=1}^\infty a_i\bigl(e_{k+1}-e_k-e_k+e_k+m(D^k)\bigr)b_i \\
&= \sum_{i=1}^\infty a_ie_{k+1}b_i=m(D_{k+1})=e_{k+1}
\end{align*}
so (\ref{eq:diag1}) for $D^{k+1}$ is satisfied.  Now if $a\in\fA_{k+1}$ then
by property (\ref{eq:diag2}) for $D_{k+1}$ we have 
$\sum_{i=1}^\infty (aa_i)\otimes b_i=
\sum_{i=1}^\infty a_i\otimes(b_ia)$, 
so it follows that $a\mult D^{k+1}= D^{k+1}\mult a$.  Now if $a\in\fA^k$, then
each $aa_i\in\fA^k$ so
\begin{align*}
a\mult D^{k+1}&=\sum_{i=1}^\infty (aa_i)\mult\bigl((e_{k+1}-e_k)\otimes(e_{k+1}-e_k)+
D^k\bigr)\mult b_i \\
&=\sum_{i=1}^\infty \bigl([aa_i(e_{k+1}-e_k)]\otimes(e_{k+1}-e_k)+(aa_i)\mult D^k
\bigr)\mult b_i \\
&=\sum_{i=1}^\infty  D^k\mult(aa_ib_i)=D^k\mult a=a\mult D^k
\end{align*}
which, by symmetric argument, is exactly the value of $D^{k+1}\mult a$.
Since any $a\in\fA^{k+1}$ is a sum $a=\pi_{k+1}(a)+(a-\pi_{k+1}(a))$
where, $\pi_{k+1}(a)\in\fA_{k+1}$ and $a-\pi_{k+1}(a)\in\fA^k$, we obtain
(\ref{eq:diag2}) for $D^{k+1}$.  \endpf

We note that to generalise our proof of the preceding result to amenable
but not contractible Banach algebras, we would require at each stage 
approximate diagonals $D^k_\alp$ such that $\norm{m(D^k_\alp)}=1$, which we
do not know how to construct, in general.  We point the reader to
\cite[Theorem 2.3]{rundes1} to see a computation performed on a Banach
algebra graded over $L_1$.

We note that we can modify the proof of Proposition 
\ref{prop:gradedlinear} to see that {\it a Banach algebra
$\fA=\ell_1\text{-}\bigoplus_{s\in F_2^1}\fA_s$ graded over $F_2^1$, 
where each $\fA_s$ is contractible with $\mathrm{AM}(\fA_s)=1$, satisfies
$\mathrm{AM}(\fA)\leq 45$}.  This is larger than
$\mathrm{AM}(F_2^1)=25$ from Example \ref{ex:flatsem}. We have found
no examples of such Banach algebras $\fA$ with
$\mathrm{AM}(\fA)>25$.  However, we conjecture only for semilattices
$S=L_n$, that {\it a Banach algebra
$\fA=\ell_1\text{-}\bigoplus_{s\in S}\fA_s$ graded over $S$, 
where each $\fA_s$ is amenable with $\mathrm{AM}(\fA_s)=1$, satisfies
$\mathrm{AM}(\fA)=\mathrm{AM}(S)$.}  It would be interesting to find
non-linear unital semilattices over which this conjecture holds.

\subsection{On allowable amenability constants}
We close by partially answering a question posed in \cite{dalesls}.  There
it is proved, that there is no
semigroup $G$ such that $1<\mathrm{AM}(G)<5$.  It is further 
conjectured that there are no semigroups $G$ for which $\mathrm{AM}(G)
\in(5,7)\cup(7,9)$.  In \cite{dalesls} there is an example given
of a noncommutative semigroup $G$ with $\mathrm{AM}(G)=7$.  For commutative
semigroups there is a further gap.

\begin{dalesq}\label{prop:dalesq}
There is no commutative semigroup $G$ such that 
\[
5<\mathrm{AM}(G)<9.
\]
\end{dalesq}

\proof Since $G$ is commutative, it is proved in \cite[Theorem 2.7]{groenbaek} 
that if $\ell^1(G)$ is amenable, then
$G$ is a Clifford semigroup, whose component groups are abelian,
graded over a finite semilattice $S$.
If $\mathrm{AM}(G)<9$, then by Theorem \ref{theo:gradedalg}
then $\mathrm{AM}(S)<9$ and hence by Theorem \ref{theo:amenconst} and
the corollary which follows it we have
\[
2|S|-1\leq \mathrm{AM}(S)\leq 5
\]
so $|S|\leq 3$.  Clearly, if $|S|=1$, $S=L_0$, and if $|S|=2$, $S=L_1$.
If $|S|=3$ then $S$ is either unital, in which case $S=L_2$, or
$S$ has 2 maximal elements, in which case $S=F_2$; in either case
$\mathrm{AM}(S)=9$, contradicting our assumptions.  Thus $S=L_0\oor L_1$.
But it then follows by a straighforward adaptation of
\cite[Theorem 2.3]{rundes1} that $\mathrm{AM}(G)
=1\oor 5$.  In particular $\mathrm{AM}(G)\leq 5$.  \endpf





\bigskip
{\bf Acknowledgements.}  The authors are grateful to H.G. Dales for
valuable questions and discussion, and the Y. Choi for providing a preprint
of his article \cite{choi}.

{
\bibliography{amenconstbib}
\bibliographystyle{plain}
}

\smallskip
Mahya Ghandehari

Address: {\sc Department of Pure Mathematics, University of Waterloo,
Waterloo, ON\quad N2L 3G1, Canada} 

E-mail: {\tt mghandeh@uwaterloo.ca}

\medskip
Hamed Hatami

Address: {\sc Department of Computer Science, University of To\-ron\-to,
Toronto, ON\quad M5S 3G4, Canada} 

E-mail:{\tt hamed@cs.toronto.ca}

\medskip
Nico Spronk

Address:  {\sc Department of Pure Mathematics, University of Waterloo,
Waterloo, ON\quad N2L 3G1, Canada}

E-mail:  {\tt nspronk@uwaterloo.ca}

\end{document}